\documentclass[11pt,amsfonts]{article}
\usepackage{graphicx}
\usepackage{latexsym}
\usepackage{amssymb}
\usepackage{amsmath}
\usepackage{enumerate}
\usepackage{dsfont}
\usepackage{xcolor}
\usepackage[urlcolor=black, colorlinks=true, linkcolor=blue, citecolor=brown]{hyperref}

\usepackage{layout}
\usepackage{eufrak}
\newtheorem{prop}{Proposition}
\newtheorem{lemma}{Lemma}
\newtheorem{definition}{Definition}

\newtheorem{theorem}{Theorem}
\newtheorem{remark}{Remark}

\def\real{{\mathord{{\rm I\kern-2.8pt R}}}}        
\def\inte{{\mathord{{\rm I\kern-2.8pt N}}}}

\def\sZZ{{\rm Z\kern-2.8ptem{}Z}}

\def\z{{\mathchoice
  {\sZZ}
  {\sZZ}
  {\rm Z\kern-0.30em{}Z}
  {\rm Z\kern-0.25em{}Z} }}
\def\sQQ{{\kern 0.27em \vrule height1.45ex width0.03em depth0em
          \kern-0.30em \rm Q}}
\def\qu{{\mathchoice
    {\sQQ}
    {\sQQ}
  {\kern 0.225em \vrule height1.05ex width0.025em depth0em \kern-0.25em \rm Q}
  {\kern 0.180em \vrule height0.78ex width0.020em depth0em \kern-0.20em \rm Q}
        }}
\def\sCC{{\kern 0.27em \vrule height1.45ex width0.03em depth0em
          \kern-0.30em \rm C}}
\def\complex{{\mathchoice
    {\sCC}
    {\sCC}
  {\kern 0.225em \vrule height1.05ex width0.025em depth0em \kern-0.25em \rm C}
  {\kern 0.180em \vrule height0.78ex width0.020em depth0em \kern-0.20em \rm C}
        }}

\newcommand{\ba}{\begin{array}}
\newcommand{\ea}{\end{array}}
\newcommand{\be}{\begin{equation}}
\newcommand{\ee}{\end{equation}}
\newcommand{\bea}{\begin{eqnarray}}
\newcommand{\eea}{\end{eqnarray}}
\newcommand{\beaa}{\begin{eqnarray*}}
\newcommand{\eeaa}{\end{eqnarray*}}

\def\z{\zeta}

\font\tenmath=msbm10 \font\sevenmath=msbm7 \font\fivemath=msbm5
\newfam\mathfam \textfont\mathfam=\tenmath
\scriptfont\mathfam=\sevenmath \scriptscriptfont\mathfam=\fivemath

\def \={{\buildrel {\rm (law)} \over =}}

\def\qed{ \hfill \vrule width.25cm height.25cm depth0cm\smallskip}

\newcommand{\basa}{\begin{assumption}}
\newcommand{\easa}{\end{assumption}}

\newcommand{\bas}{\begin{assum}}
\newcommand{\eas}{\end{assum}}

\newcommand{\ignore}[1]{}
\begin{document}

\renewcommand{\thefootnote}{\fnsymbol{footnote}}

\renewcommand{\thefootnote}{\fnsymbol{footnote}}

\title{Multifractal random walks with fractional Brownian motion via Malliavin calculus }
\author{ Alexis Fauth $^{1}$ \footnote{Head of quantitative research at Invivoo, 13 rue de l'Abreuvoir, 92400 Courbevoie, France. }\qquad %
Ciprian A. Tudor $^{2,3}$ \footnote{Supported by the CNCS grant PN-II-ID-PCCE-2011-2-0015. Associate member of the team Samm, Universit\'e de Panth\'eon-Sorbonne Paris 1 }\vspace*{0.1in} \\
$^{1}$SAMM, Universit\'e de Paris 1 Panth\'eon-Sorbonne,\\
90, rue de Tolbiac, 75634, Paris, France. \\
 \quad alexis.fauth@invivoo.com \vspace*{0.1in} \\
 $^{2}$ Laboratoire Paul Painlev\'e, Universit\'e de Lille 1\\
 F-59655 Villeneuve d'Ascq, France.\\
 $^{3}$ Department of Mathematics, \\
Academy of Economical Studies, Bucharest, Romania \vspace*{0.1in}\\
 \quad tudor@math.univ-lille1.fr\vspace*{0.1in}}

\maketitle

\begin{abstract}
We introduce a Multifractal Random Walk (MRW) defined as a stochastic integral of an infinitely divisible noise with respect to a dependent fractional Brownian motion. Using the techniques of the Malliavin calculus, we study the existence of this object and its properties. We then propose a continuous time model in finance that captures the main properties observed in the empirical data, including the leverage effect. We illustrate our result by numerical simulations.
\end{abstract}

\vskip0.5cm

{\bf  2010 AMS Classification Numbers:} 60C30, 60H07, 60H05.

 \vskip0.3cm


\section{Introduction}

Starting with the seminal work of Mandelbrot about cotton price \cite{M}, several studies of financial stock prices times series, have allowed to exhibit some particularities of their fluctuations. Without making a comprehensive list, we can mention the appearance in empirical data  of the following properties:  non Gaussian distributions due to now well known fat tails  of financial returns, the so-called volatility clustering that means that the volatility fluctuations are of intermittent and correlated nature, scaling invariance, long run correlation in volatility, leverage effect and so on (see e.g \cite{GPVAMS}, \cite{CMPF}, \cite{Voit} for an extensive review). Thus, constructing theoretical models for  financial returns that include all the properties listed before  appears as a very interesting challenge. Many scientific works, in economics or mathematics, proposed various  models for asset returns.  The  ARCH model introduced by Engle \cite{Engle}  offers an interesting base of work, and after the seminar work by Engle a vast literature on ARCH and related models has been developed. One of the  first extensions of the ARCH model, called the GARCH model,  has been introduced by Bollerslev in \cite{Bol}  and it also has been the object of various generalizations. While the GARCH models are capable of capturing volatility clustering, there are some drawbacks of the model. For example, the GARCH models  are unable to represent volatility asymmetry.  Due to the presence of the squared observed data in the conditional variance equation, the positive and negative values of the lagged innovations have the same effect on the conditional variance. In the finance literature, it has been noticed that volatility often is affected by negative and positive shocks in different ways. Another inconvenient of the  GARCH model is the fact that, to ensure positiveness of the conditional variance, non-negative constraints on the  coefficients in the variance equation must be imposed. To take into account the asymmetric effects on conditional second moments and to avoid the non-negativity constraints on the coefficients in the variance equation, in \cite{Nelson} the so-called EGARCH model has been proposed by D. Nelson.

On the other hand, many recent empirical studies, based on huge data sets available nowadays, put in light new aspects of the financial returns.  For example, they  suggested that the fluctuations of the asset process displays multifractal properties (see e.g. \cite{FCM} or \cite{GBPTD}). Taking into account these multifractal character, several authors proposed models based on the on "cascade" random processes and  "multifractal random walk (MRW)". We refer, among others, to \cite{CF}, \cite{BMD} or \cite{BKM}.
Usually, the noise in these models are defined by
\begin{equation}\label{1}
Z(t)= Y(X([0,t]), \hskip0.5cm t\in \mathbb{R}
\end{equation}
where $X$ is a multifractal random measure and $Y$ is a self-similar process with stationary increments, independent of $X$. Another construction defines the noise via a stochastic integral as (see e.g. \cite{BaMa}, see \cite{ACCP})
\begin{equation}
\label{2}
Z(t)= \int_{0} ^{t} Q(u) \mathrm{d}Y(u)
\end{equation}
where $Q$ is a suitable fractal noise and  $Y$ is a self-similar process with stationary increments, independent of $Q$. A natural choice for the process $Y$, both in (\ref{1}) and (\ref{2}),  is the fractional Brownian motion (fBm).

While  these multifractal models (\ref{1}), (\ref{2}) are able to capture several properties observed on the empirical data (scaling, volatility clustering or long-range dependance), they reproduce poorly the leverage effect. The leverage effect is understood as the correlation between the log-return at a fixed time $t$ and his volatility (that may be defined as the squared or absolute log-return) varying around $t$, for example $t-100\Delta t\leq t'\leq t+100\Delta t$. Empirical studies (see e.g. \cite{Black} or \cite{BMP}) have shown that this quantity is close to 0 for past volatility and follow an exponential law for future volatility. This stylized fact is interpreted as a panic effect, the two quantities are negatively correlated just after $t$ and go back quickly to 0 as $t'$ increases.

 We introduce here a generalization of the model (\ref{2}), able to take into account the leverage effect and allows a more  flexibility for long-range dependence in log-return. That is, in the model (\ref{2}) we will consider the processes $Q$ to be multifractal process called "infinitely divisible cascading noise" and the process $Y$ to be a fractional Brownian motion with Hurst parameter $H>\frac{1}{2}$.  The same situation has been treated in \cite{ACCP}, \cite{Lud}, but in addition, we construct a model in which the processes $Q$ and $B^{H}$  are not independent anymore. An alternative construction has been proposed in \cite{BDM} but in this reference the integral in (\ref{2}) is not a true stochastic integral with respect to the fBm. We will define the integral (\ref{2}) using the techniques of the Malliavin calculus and we will study the properties of the fractional Multifractal Random Walk.

 Our paper is organized as follows. Section 2 contains some preliminaries on multifractal process and the basic elements of the Malliavin calculus that we will need in our paper. In Section 3 we analyze the stochastic integrals that are used to define the MRW while in Section 4 we study the existence of the fractional MRW as a limit of a family of stochastic integrals. The properties of the fractional MRW (scaling, moments etc) are discussed in Section 5. The last section (Section 6) contains a numerical analysis of the data for several financial indices and we compare the simulation of our theoretical model with the real data.

\section{Preliminaries}

We present here he basic facts related to the infinitely divisible cascading noises and we introduce the basic tools  of the Malliavin calculus.

\subsection{Infinitely divisible cascading noise}

Let $M$ denote an infinitely divisible, independently scattered random measure on the set $\mathbb{R} \times \mathbb{R}_{+}$  with generating infinitely divisible distribution $G $ satisfying
\begin{equation*}
\int_{\mathbb{R}} e^{qx}G(\mathrm{d}x)= e^{-\rho (q)}
\end{equation*}
for some function $\rho$ and for every $q\in \mathbb{R}$. We assume that $M$ has  control measure $m$ on $\mathbb{R} \times \mathbb{R} _{+}$ meaning that for every Borel set $A \subset \mathbb{R} \times \mathbb{R} _{+}$ it holds
\begin{equation*}
\mathbb{E} e^{qM(A) } =e^{-\rho (q) m(A)} \mbox{ for every } q\in \mathbb{R}.
\end{equation*}
The fact that $M$ is independently scattered  means that the random variables
$$M(A_{1}), M(A_{2}),..., M(A_{n}) $$
are independent whenever the Borel sets $A_{1},..., A_{n} \in \mathbb{R} \times \mathbb{R} _{+}$ are disjoint.
We define the {\it Infinitely Divisible Cascading noise} (IRC)  by
\begin{equation}\label{irc}
Q_{r}(t)= \frac{ e^{M (C_{r}(t))}}{ \mathbb{E}e^{M (C_{r}(t))}}
\end{equation}
for every  $r>0$ and $t\in \mathbb{R}$. Here $C_{r}(t)$ is the cone in $\mathbb{R} \times \mathbb{R} _{+}$ defined by
\begin{equation}
\label{crt}
C_{r}(t)=\{ (t',r'), r\leq r' \leq 1, t-\frac{r'}{2} \leq t' \leq t+ \frac{r'}{2} \}.
\end{equation}

We will use the  following facts throughout our paper. We refer to \cite{CRA} or \cite{ACCP} for the their proofs. First, let us note the scaling property of the moments of the IRC
\begin{equation*}
\mathbb{E} Q_{r}(t) ^{q} = e^{-\varphi(q) m_{r}(0) }
\end{equation*}
and the expression of its   covariance: for every $r>0$ and $t,s\in \mathbb{R}$
\begin{equation}\label{covQ}
\mathbb{E} Q_{r}(t) Q_{r}(s)= e^{-\varphi (2) m_{r} (\vert t-s\vert ) }
\end{equation}
where we denoted by for $u\geq 0, r>0$.
\begin{equation*}
m_{r}(u) =m\left( C_{r}(0) \cap C_{r}(u) \right)
\end{equation*}
and by
\begin{equation}\label{fi}
\varphi (q)= \rho(q)-q\rho (1).
\end{equation}

The scaling of the moment of $Q$ can be extended to the following scaling property in distribution : for $t\in (0,1)$
\begin{equation}
\label{sq}
\left( Q_{rt} (tu) \right) _{u\in \mathbb{R}} =^{(d)} e^{\Omega _{t} }\left( Q_{r}(u) \right) _{u\in \mathbb{R}}
\end{equation}
where $"=^{(d)}"$ means equivalence of finite dimensional distributions. Here $\Omega _{t}$ denotes a random variable independent by $Q$, which satisfies, if the measure $m$ is given by (\ref{m}),
\begin{equation}
\label{o}
\mathbb{E} e^{q\Omega _{t}} = t^{q \varphi (q)}.
\end{equation}

\begin{remark}
As noticed in \cite{ACCP}, we have $\varphi (2)<0.$
\end{remark}

In \cite{ACCP} (see also \cite{Lud} ) the {\it Multifractal random walk } (MRW) based on fractional Brownian motion is defined as limit when $r\to 0$ (in some sense) of the family of  stochastic integrals $\left( Z^{H}_{r}(t)\right) _{r>0} $ defined by
\begin{equation}\label{zrt}
Z^{H} _{r} (t) = \int_{0} ^{t} Q_{r}(u) \mathrm{d}B^{H}(u), \hskip0.5cm t\in [0,T]
\end{equation}
where $(B^{H}_{t})_{t\in [0,T]}$ is a fractional Brownian motion with Hurst parameter $H\in (0,1)$. The fractional Brownian motion $(B^{H}_{t}) _{t\in [0,T]}$ with  Hurst parameter $H\in (0, 1)$ is a centered Gaussian process starting from zero with covariance function $$R^{H}(t,s):= \frac{1}{2} \left( t^{2H}+s^{2H} -\vert t-s\vert ^{2H}\right), \hskip0.5cm s,t \in [0,T].$$ In \cite{ACCP} it is assumed that $M$ and $B^{H}$ {\bf are independent. } Therefore the stochastic integral with respect to $B^{H}$ in (\ref{zrt}) behaves mainly as a Wiener integral since, because of the independence, the integrand $Q_{r}(u)$ can be viewed as deterministic function for the integrator $B^{H}$.

Another important fact in the development of this theory is that the IRC $Q$ is a martingale with respect to the argument $r$. Let us recall the following result (see \cite{CRA}):
\begin{lemma}
For every $u>0$ the stochastic process $\left( Q_{r}(u) \right) _{r>0}$ is a martingale with respect to its own filtration. As a consequence,  for every $u,v,r,r'>0$ with $r<r'$ it holds
\begin{equation}
\label{mart}
\mathbb{E}Q_{r}(u) Q_{r'}(v)= \mathbb{E}Q_{r}(u) Q_{r}(v).
\end{equation}
\end{lemma}
The property (\ref{mart}) plays an important role in the construction of the MRW process in \cite{ACCP} or \cite{CRA}.

\subsection{Malliavin calculus}
 Let $(W_{t})_{t\in T}$ be a classical
Wiener process on a standard Wiener space $\left( \Omega
,{\mathcal{F}},\mathbf{P}\right) $. By $W(\varphi)$ we denote the Wiener integral of the function $\varphi \in L^{2}(T)$ with respect to the Brownian motion $W$.
We denote by $D$  the Malliavin  derivative operator that acts on smooth functionals of the form $F=g(W(\varphi _{1}), \ldots , W(\varphi_{n}))$ (here $g$ is a smooth function with compact support and $\varphi_{i} \in L^{2}(T)$ for $i=1,..,n$)
\begin{equation*}
DF=\sum_{i=1}^{n}\frac{\partial g}{\partial x_{i}}(W(\varphi _{1}), \ldots , W(\varphi_{n}))\varphi_{i}.
\end{equation*}
The operator $D$ can be extended to the closure $\mathbb{D}^{p,2}$ of smooth functionals with respect to the norm
\begin{equation*}
\Vert F\Vert _{p,2}^{2} = \mathbb{E}F^{2}+ \sum_{i=1}^{p} \mathbb{E} \Vert D^{i} F\Vert ^{2} _{L^{2}T i}
\end{equation*}
where the $i$ th Malliavin derivative $D^{(i)}$ is defined iteratively.  The adjoint of
$D$ is denoted by $\delta $ and is called the divergence (or
Skorohod) integral. Its domain ($Dom(\delta)$ ) coincides with the class of stochastic processes $u\in L^{2}(\Omega \times T)$ such that
\begin{equation*}
\left| \mathbb{E}\langle DF, u\rangle \right| \leq c\Vert F\Vert _{2}
\end{equation*}
for all $F\in \mathbb{D}^{1,2}$ and $\delta (u)$ is the element of $L^{2}(\Omega)$ characterized by the duality relationship
\begin{equation*}
\mathbb{E}(F\delta (u))= \mathbb{E}\langle DF, u\rangle.
\end{equation*}
For adapted integrands, the divergence integral coincides to
the classical It\^o integral. A subset of $Dom (\delta)$ is the space $\mathbb{L} ^{1,p}$ of the stochastic processes such that $u_{t}$ is Malliavin differentiable for every $t$ and
\begin{equation*}
\Vert u\Vert _{1,p}^{p}:= \Vert u\Vert ^{p} _{L^{p}(T\times \Omega) }+ \Vert Du \Vert ^{p} _{L^{p}(T ^{2} \times \Omega)}<\infty.
\end{equation*}
We will need Meyer's inequality that allows to estimate the $L^{p} $ moment of the Skorohod integral
\begin{equation}
\label{meyer}
\mathbb{E} \left| \delta (u) \right| ^{p} \leq \Vert u\Vert _{1,p} ^{p}.
\end{equation}
For the $L^{2}$ moment of the Skorohod integral we have the explicit formula
\begin{equation}
\label{l2}
\mathbb{E} \delta (u)^{2} = \int_{T} u_{s} ^{2} \mathrm{d}s + \int _{T} \int_{T} D_{r} u_{s} D_{s} u_{r}\mathrm{d}r\mathrm{d}s
\end{equation}
if $u\in \mathbb{L} ^{1,2}$. We also recall that the Malliavin derivative satisfies  the chain rule
\begin{equation}
\label{chain}
Df(F)= f'(F)DF
\end{equation}
if $f$ is a differentiable function and $F\in \mathbb{D}^{1,2}$.

\section{The construction of the fractional Multifractal Random Walk with dependent noise}
Our purpose is to give a meaning to the stochastic integral (\ref{zrt}) in the situation when the IRC $Q$ and the fBm $B^{H}$ are not independent. We will use techniques related to the Malliavin calculus. In order to apply these type of techniques we will restrict ourselves to the case when the measure $M$ introduced in Section 2 is Gaussian.

\subsection{The Gaussian isonormal noise}
We introduce a MRW  without independence between the measure  $M$ (denote by $W$ in our settings) and the integrator $B^{H}$. We will restrict to the case where $M$ is a Gaussian measure. More precisely, we will consider $(W(h), h\in H)$ an isonormal process, that is, a centered Gaussian family with
$$\mathbb{E} W(h)  W(g) = \langle h,g \rangle _{H}$$
for every $g,h \in H$. The Hilbert space $H$ will be
$$H= L^{2} \left( \mathbb{R}\times \mathbb{R} _{+} , {\cal{B}} (\mathbb{R} \times \mathbb{R} _{+}) , m)\right)$$
where
$m$ is the control measure. In this work we will consider
\begin{equation}\label{m}
m(dt, dr)= dt \frac{c dr}{r^{2}} \mbox{ if } 0<r\leq 1
\end{equation}
and $m$ vanishes if $r\geq 1$. Here $c$ is a strictly positive constant. This is called in \cite{BM1} (see also \cite{ACCP}) the exact invariant scaling case.

  The following properties of the noise $W$ are immediate: denote $W(A)= W(1_{A}) $ for $A \in {\cal{B}} (\mathbb{R} \times \mathbb{R} _{+}) $.
  \begin{description}
    \item{$\bullet $}
    For every $A \in {\cal{B}} (\mathbb{R} \times \mathbb{R} _{+}) $ it holds that
    \begin{equation*}
    W(A) \sim N(0, m(A) ).
    \end{equation*}

  \item{$\bullet $}We have
  \begin{equation*}
  \mathbb{E} W(A) W(B) = 0
  \end{equation*}
  if the Borel sets $A,B \subset \mathbb{R} \times \mathbb{R} _{+}$ are disjoint. This implies that the random variables $W(A) $ and $W(B)$ are independent, so the random measure $W$ is independently scattered.

  \item{$\bullet $} For every $q\in \mathbb{R}$ and $A \in {\cal{B}} (\mathbb{R} \times \mathbb{R} _{+}) $
  we have
  \begin{equation*}
  \mathbb{E} e^{qW(A) }= e^{-\frac{1}{2} q^{2} m(A)   }
  \end{equation*}
  which means that
  \begin{equation*}
  \rho (q)= -\frac{1}{2} q^{2} \mbox{ and } \varphi (q)= -\frac{1}{2} q^{2} + \frac{q}{2} =-\frac{1}{2}q(q-1).
  \end{equation*}
  \end{description}

  Let us use the following notation:
  $$m_{1}(dt) = dt \mbox{ and } m_{2}(dr) = c\frac{dr}{r^{2}}.$$
  Above $m_{1}$ is the Lebesque measure on $\mathbb{R}$ and $m_{2}$ is a measure on $\mathbb{R}_{+}$. Clearly $m= m_{1} \otimes m_{2}$, the product measure.

  For every $t\geq 0$ and $ A_{2} \in {\cal{B}} (\mathbb{R} _{+}) $ such that
  $$m_{2} (A_{2}) =1$$ (take  for example $A_{2}= (c, \infty)$) 
we set

  \begin{equation}\label{w1}
W^{(1) }(t): =   W (1_{[0,t] \times A_{2}} ).
  \end{equation}

  The following result is immediate.
  \begin{prop}\label{p1}
  If $A_{2} \in {\cal{B}} (\mathbb{R}_{+} ) $ is such that $m_{2}(A_{2})=1$ then the process $(W^{(1)}(t) ) _{t\geq 0}$ given by (\ref{w1}) is a standard Brownian motion on the same probability space as $W$.
  \end{prop}
  {\bf Proof: } It is clear that $W^{(1)}$ is a Gaussian process. Let us compute its covariance. For every $s,t \geq 0$ it holds that
  \begin{eqnarray*}
  \mathbb{E} W^{(1)}(t) W^{(1)}(s) &=& \mathbb{E} W (1_{[0,t] \times A_{2}} )W (1_{[0,s] \times A_{2}} ) \\
  &=&\langle 1_{[0,t] \times A_{2}},1_{[0,s] \times A_{2}}\rangle _{H} =\langle 1_{[0,t]}, 1_{[0,s]}\rangle _{L^{2}(\mathbb{R})}m_{2}(A_{2})\\
  &=& t\wedge s
  \end{eqnarray*}
and this implies that $W^{(1)}$ is a Wiener process with respect to its own filtration. \qed

\vskip0.3cm

\subsection{The approximating  Multifractal Random Walk with dependent fractional Brownian motion}

We introduce the Multifractal Random Walk with based on the fractional Brownian motion by  the formula
\begin{equation}\label{mrw-2}
Z^{H} _{r} (t)= \int_{0} ^{t} Q_{r} (u) \mathrm{d}B^{H} (u)
\end{equation}
where for every $r>0, u\geq 0$ the integrands $Q_{r}(u)$ is defined by (\ref{irc}) with $W$ instead of $M$ and the fractional Brownian motion $B^{H}$ is given by
\begin{equation}
B^{H}(t)= \int_{0} ^{t} K^{H}(t,s) \mathrm{d}W^{(1)} (s)
\end{equation}
where
$K^{H}$ is the usual kernel of the fractional Brownian motion and $W^{(1)}$ is the Brownian motion defined by (\ref{w1}). By Proposition \ref{p1}, it is clear that $B^{H}$ is a fractional Brownian motion.  Recall that, when $H>\frac{1}{2}$ the kernel $
K^{H}\left( t,s\right) $ has the expression
$$K^{H}(t,s)=c_{H}s^{1/2-H}%
\int_{s}^{t}(u-s)^{H-3/2}u^{H-1/2}\mathrm{d}u$$
(see \cite{N}) where $t>s$ and $c_{H}=\left( \frac{%
H(2H-1)}{\beta (2-2H,H-1/2)}\right) ^{1/2}$ and $\beta (\cdot ,\cdot )$ is
the Beta function. For $t>s$, the kernel's derivative is
$$\frac{\partial
K^{H}}{\partial t}(t,s)=c_{H}\left( \frac{s}{t}\right) ^{1/2-H}(t-s)^{H-3/2}.$$
In the sequel we will simply denote $K^{H}:=K$.

The stochastic integral in (\ref{mrw-2}) is a divergence (Skorohod) integral with respect to $B^{H}$ as defined in e.g. \cite{N}. Actually, when $H>\frac{1}{2}$ we can write
$$Z^{H}_{r}(t) = \int_{0} ^{t} Q_{r}(s) \mathrm{d}B^H(s)= \int_{0} ^{t} dW^{(1) }(s)\left( \int_{s}^{t} \mathrm{d}a \partial _{1} K(a,s)  Q_{r} (a)\right).$$
The last equality is due to the definition of the stochastic integral with respect to $B^{H} $ (see \cite{N}, Chapter 5 for example). We can also  express  $Z_{r}(t)$ as a Skorohod  integral with respect to the measure $W$ by the formula
\begin{equation}\label{zrt2}
Z^{H}_{r}(t) = \int_{0} ^{t} \int_{A_{2}} \mathrm{d}W (s,r')\left(  \int_{s}^{t} \mathrm{d}a \partial _{1} K(a,s)  Q_{r} (a)\right)
\end{equation}
where $A_{2}$ is the Borel set satisfying $m_{2}(A_{2})=1$ and that appears in the definition of $W^{(1)}$.

\begin{remark}
 Clearly the fBm $B^{H}$ and the noise $W$ are dependent. Actually, it can be seen that  $B^{H}$ and $W$ are correlated in general. Indeed, for every Borel set $A \subset \mathbb{R} \times \mathbb{R}_{+}$ and for every $t\geq 0$
\begin{eqnarray*}
\mathbb{E} B^{H}(t) W(A)&=&  \mathbb{E} \int_{0} ^{t} K(t,s) \mathrm{d}W^{(1)}(t) W(A) \\
&=& \mathbb{E} \int_{0} ^{t} \int_{0} ^{1} K(t,s) \mathrm{d}W(s,r) W(A) \\
&=& \mathbb{E} \int \int _{ \left( (0,t) \times (1, \infty ) \right) \cap A } K(t,s) \mathrm{d}s \frac{\mathrm{d}r}{r^{2} }
\end{eqnarray*}
and this is in general not zero.

\end{remark}

Let us further discuss the dependence between the integrator $B^{H}$ and the integrand $Q$ in (\ref{zrt}). We need to distinguish two situations.  It depends on the relations of the set $A_{2}$ and the unit interval $(0,1)$.

\vskip0.3cm

{\bf The disjoint case: } We have $A_{2} \cap (0,1)= \emptyset$ (this happens when $A_{2}= (c, \infty)$ with $c>1$ for example).

On the other hand, $B^{H}_{s}$ is independent by $W(C_{r} (t))$ for every $s,t, r$.  Indeed, since $C_{r} (t) \subset \mathbb{R} \times (0,1) $ (see (\ref{crt})) we have
$$\mathbb{E}B^{H}_{s} W(C_{r}(t))=0$$
and since $(B^{H}(s), W (C _{r}(t))) $ is a Gaussian vector, we obtain the independence.

\vskip0.3cm

 {\bf The non-disjoint case: } $A_{2} \cap (0,1)\not= \emptyset$ (this happens when $A_{2}= (c, \infty)$ with $c<1$). In this case $\mathbb{E}B^{H}(s) W(C_{r}(t))\not=0$ in general and so $B^{H}(s)$ and $W(C_{r}(t))$ are dependent.

 \vskip0.2cm

 We will refer throughout this work to the two situations above as the {\it the disjoint case} and the {\it non-disjoint case. } Basically, the results in the disjoint case can be obtained by following the arguments in \cite{ACCP} while in the non-disjoint case the context is different because of the appearance of the Malliavin derivatives in the expression of square mean of (\ref{zrt}).

\begin{prop}\label{p2}Suppose $H>\frac{1}{2}$.
Then for every $r>0$ the stochastic Skorohod integral in (\ref{mrw-2}) is well-defined.

\end{prop}
{\bf Proof: } We will use the representation of (\ref{zrt}) of  $Z^{H}_{r}(t)$ as a Skorohod  integral with respect to the measure $W$.
Next we use the bound (\ref{meyer}) with $p=2$. Let us apply it to the process
 $$(s,r') \to 1_{[0,t] }(s)1_{A_{2}} (r)\int_{s}^{t} \mathrm{d}a \partial _{1} K(a,s)  Q_{r} (a)$$
(which is a two-parameter process) and to $H= L^{2} (\mathbb{R}\times \mathbb{R}_{+}, m)$. The variables $t,r$ and the set $A_{2}$ with $m_{2}(A_{2})=1$ are fixed. Below $D$ is the Malliavin derivative with respect to the isonormal process $W$ (see Section 2). We will get, by (\ref{meyer}),
\begin{eqnarray*}
\mathbb{E}( Z^{H} _{r} (t) ) ^{2} &\leq & \mathbb{E} \int_{0} ^{t} \mathrm{d}s \int_{A_{2}} \frac{c \mathrm{d}r'}{(r') ^{2}} \left( \int_{s}^{t} \mathrm{d}a \partial _{1} K(a,s)  Q_{r} (a) \right) ^{2} \\
&&+ \mathbb{E} \int_{0} ^{t} \mathrm{d}s \int_{A_{2}} \frac{c \mathrm{d}r'}{(r') ^{2}}\int_{0} ^{t} \mathrm{d}\alpha  \int_{A_{2}} \frac{ c \mathrm{d}\beta} {\beta ^{2}} \left( D_{\alpha , \beta } \int_{s}^{t} \mathrm{d}a \partial _{1} K(a,s)  Q_{r} (a) \right)^{2} \\
&=&\mathbb{E} \int_{0} ^{t} \mathrm{d}s \int_{A_{2}} \frac{c \mathrm{d}r'}{(r') ^{2}} \left( \int_{s}^{t} \mathrm{d}a \partial _{1} K(a,s)  Q_{r} (a) \right) ^{2}\\
&&+ \mathbb{E} \int_{0} ^{t} \mathrm{d}s \int_{A_{2}} \frac{c \mathrm{d}r'}{(r') ^{2}}\int_{0} ^{t} \mathrm{d}\alpha  \int_{A_{2}} \frac{ c \mathrm{d}\beta} {\beta ^{2}} \left( \int_{s}^{t} \mathrm{d}a \partial _{1} K(a,s)  D_{\alpha , \beta }Q_{r} (a) \right)^{2}\\
&=&\mathbb{E} \int_{0} ^{t} \mathrm{d}s \int_{A_{2}} \frac{c \mathrm{d}r'}{(r') ^{2}} \left( \int_{s}^{t} \mathrm{d}a \partial _{1} K(a,s)  Q_{r} (a) \right) ^{2}\\
&&+\mathbb{E} \int_{0} ^{t} \mathrm{d}s \int_{A_{2}} \frac{c \mathrm{d}r'}{(r') ^{2}}\int_{0} ^{t} \mathrm{d}\alpha  \int_{A_{2}} \frac{ c \mathrm{d}\beta} {\beta ^{2}} \left( \int_{s}^{t} \mathrm{d}a \partial _{1} K(a,s)  Q_{r} (a)1_{C_{r}(a)}(\alpha , \beta)  \right)^{2}\\
&:=& T_{1}+ T_{2}.
\end{eqnarray*}
since, by the chain rule of the Malliavin operator (\ref{chain})
\begin{equation}
\label{dq}
 D_{\alpha , \beta }Q_{r} (a) = Q_{r}(a) 1_{C_{r}(a) } (\alpha , \beta ).
 \end{equation}

We need again to consider two cases.

\vskip0.2cm

{\bf The disjoint case: } We have $A_{2} \cap (0,1)= \emptyset$ (this happens when $A_{2}= (c, \infty), c>1$ for example).
In this case the term denoted by $T_{2} $ vanishes because
$$1_{C_{r}(a) } (\alpha , \beta ) 1_{A_{2}}(\beta )= 0.$$
We need to show that $$T_{1}=m_{2}(A_{2})\mathbb{E} \int_{0} ^{t} \mathrm{d}s  \left( \int_{s}^{t} \mathrm{d}a \partial _{1} K(a,s)  Q_{r} (a) \right) ^{2}$$  is finite under condition (\ref{cc}) and this is exactly the computation in \cite{ACCP}, proof of Proposition 2.1.\\

{\bf The non-disjoint case: }Suppose  $A_{2} \cap (0,1)\not= \emptyset$ (this happens when $A_{2}= (c, \infty)$ with  $c<1$).
In this case we need also to show that
$$T_{2} <\infty$$
By bounding the indicator function $1_{C_{r}(a)}(\alpha , \beta)$ by 1 (note that $\partial_{1} K(a,s) \geq 0$ for every $a,s$)
\begin{eqnarray*}
&&\mathbb{E}\left( \int_{s}^{t} \mathrm{d}a \partial _{1} K(a,s)  Q_{r} (a)1_{C_{r}(a)}(\alpha , \beta)  \right)^{2}\\
&\leq & \mathbb{E}\left( \int_{s}^{t} \mathrm{d}a \partial _{1} K(a,s)  Q_{r} (a)  \right)^{2}
\end{eqnarray*}
and by computing the integrals $dr'$ and $d\beta$
\begin{eqnarray*}
T_{2} &\leq &m_{2} (A_{2}) ^{2}\mathbb{E}\int_{0} ^{t} \mathrm{d}s \int_{0} ^{t} d\alpha \left( \int_{s}^{t} \mathrm{d}a \partial _{1} K(a,s)  Q_{r} (a) \right)^{2} \\
&=& m_{2}(A_{2}) ^{2}\mathbb{E}\int_{0} ^{t} \mathrm{d}s \left( \int_{s}^{t} \mathrm{d}a \partial _{1} K(a,s)  Q_{r} (a) \right)^{2}\\
&=&  t T_{1}
\end{eqnarray*}
and this is finite under (\ref{cc}).
\qed


\section{The Multifractal Random Walk}

The purpose of this section is to study the limit as $r\to 0$ of the family of stochastic processes $\left(Z^{H}_{r}(t)\right)$ with fixed $t>0.$ This limit will be called the {\it the Multifractal Random Walk.  }

We will assume throughout this paragraph that $Z^{H}_{r}(t)$ is defined by (\ref{mrw-2}) with $W^{(1)}$ given by (\ref{w1}). Moreover we will suppose
$$m_{2}(A_{2}) = 1$$
and $A_{2} \cap (0,1) \not= \emptyset.$ The  disjoint case  $A_{2} \cap (0,1) = \emptyset$ follows from \cite{ACCP}.

We will consider the following assumption:

\begin{equation}
\label{cc}
c\varphi (2)+2H >1.
\end{equation}

\begin{remark} \label{rem1}
Since  $\varphi (2)= -1$, the condition (\ref{cc}) means that $c<2H-1$. Since $H>\frac{1}{2}$ we can choose $c>0$ in order to have (\ref{cc}). By assuming (\ref{cc}), we are in the case (A) in \cite{ACCP}.
\end{remark}

We have the following limit theorem.

\begin{theorem}
\label{t1}
Assume (\ref{cc}). For every $t>0$ the sequence  of stochastic process $\left (Z^{H}_{r}(t)\right) _{r>0}$ defined by (\ref{mrw-2}) converges in $L^{2}(\Omega)$ to a random variable $Z^{H}(t)$.
\end{theorem}
{\bf Proof: } Let us fix $r, r' \in (0,1)$ with $r<r'.$ We will first compute the $L^{2}(\Omega)$ norm of the increment $Z^{H}_{r}(t)-Z^{H}_{r'}(t)$ where $t>0$ is fixed. We can write, with $W^{(1)}$ given by (\ref{w1})
\begin{eqnarray*}
\mathbb{E} \left| Z^{H}_{r}(t)-Z^{H}_{r'}(t)\right| ^{2}
&=& \mathbb{E} \left[ \int_{0} ^{t} \left( Q_{r} (u) -Q_{r'}(u) \right) \mathrm{d}B^{H} _{u} \right] ^{2}  \\
&=& \mathbb{E} \left[ \int_{0} ^{t} \left( \int_{u}^{t} \left( Q_{r}(a)-Q_{r'}(a) \right) \partial _{1} K(a,u) \right) \mathrm{d}W^{(1)} _{u} \right] ^{2} \\
&=&  \mathbb{E} \left[ \int_{0} ^{t} \int_{A_{2}}\left( \int_{u}^{t}\mathrm{d}a \left( Q_{r}(a)-Q_{r'}(a) \right) \partial _{1} K(a,u) \right) \mathrm{d}W(u,x) \right] ^{2} \\
&\leq & \mathbb{E} \int_{0} ^{t} \mathrm{d}u \int_{A_{2}} \frac{c \mathrm{d}x}{x^{2}}  \left( \int_{u}^{t}\mathrm{d}a \left( Q_{r}(a)-Q_{r'}(a) \right) \partial _{1} K(a,u) \right)^{2} \\
&&+ \mathbb{E}  \int_{0} ^{t} \mathrm{d}u \int_{A_{2}} \frac{c \mathrm{d}x}{x^{2}} \int_{0} ^{t} \mathrm{d}v \int_{A_{2}} \frac{c \mathrm{d}y}{y^{2}}
\left[ D_{v,y}  \int_{u}^{t}\mathrm{d}a \left( Q_{r}(a)-Q_{r'}(a) \right) \partial _{1} K(a,u) \right] ^{2}
\end{eqnarray*}
where  we used the bound (\ref{meyer}) for the $L^{2}$ norm of the divergence operator. Using the differentiation rule  (\ref{dq}), we get
\begin{eqnarray*}
&&\mathbb{E} \left| Z^{H}_{r}(t)-Z^{H}_{r'}(t)\right| ^{2} \\
&\leq & m_{2} (A_{2}) \int_{0} ^{t} du \int_{u}^{t} \mathrm{d}a \int_{u}^{t} \mathrm{d}b \partial _{1} K(a,u)\partial _{1} K(b,u)\\
&&\times \mathbb{E} \left( Q_{r}(a) -Q_{r'}(a) \right)  \left( Q_{r}(b) -Q_{r'}(b) \right)\\
&&+ m_{2}(A_{2}) \int_{0} ^{t} \mathrm{d}u \int_{0} ^{t} \mathrm{d}v \int_{A_{2}} \frac{c \mathrm{d}y}{y^{2}}  \left[ \int_{u}^{t} \mathrm{d}a\partial_{1}K(a,u) \left( Q_{r}(a) 1_{C_{r}(a)} (v,y) - Q_{r'}(a) 1_{C_{r'}(a)} (v,y)\right) \right]^{2}\\
&:=& A(r,r')+ B(r,r').
\end{eqnarray*}
Let us first treat the term denoted by $A(r,r').$ Using property (\ref{mart}),
\begin{eqnarray*}
A(r,r')&=&m_{2} (A_{2}) \int_{0} ^{t} \mathrm{d}u \int_{u}^{t} \mathrm{d}a \int_{u}^{t} \mathrm{d}b \partial _{1} K(a,u)\partial _{1} K(b,u)\\
&&\times \mathbb{E} \left[ Q_{r}(a) Q_{r}(b)-2Q_{r}(a) Q_{r}(b)+ Q_{r'}(a) Q_{r'} (b) \right] \\
&=& m_{2} (A_{2}) \int_{0} ^{t} du \int_{u}^{t} \mathrm{d}a \int_{u}^{t} \mathrm{d}b \partial _{1} K(a,u)\partial _{1} K(b,u)\mathbb{E}\left[ Q_{r'}(a) Q_{r'} (b) -Q_{r}(a) Q_{r} (b)  \right]\\
&=& \mathbb{E} \int_{0} ^{t} \mathrm{d}u \left( \int_{u}^{t} \mathrm{d}a \partial _{1} K (a,u) Q_{r'} (a) \right) ^{2} \mathrm{d}u -\mathbb{E} \int_{0} ^{t} \mathrm{d}u \left( \int_{u}^{t} \mathrm{d}a \partial _{1} K (a,u) Q_{r} (a) \right) ^{2} \mathrm{d}u\\
&=& \mathbb{E} \left( \int_{0} ^{t} Q_{r'}(u) \mathrm{d}\tilde{B} ^{H}_{u} \right) ^{2}-\mathbb{E} \left( \int_{0} ^{t} Q_{r}(u) \mathrm{d}\tilde{B} ^{H}_{u} \right) ^{2}
\end{eqnarray*}
where $\tilde{B} ^{H}$ denotes a fractional Brownian motion independent by $W$.  The convergence of this term is similar to the study in \cite{ACCP} we are in the case (A) in their paper (see Remark \ref{rem1}).

The summand $B(r,r')$ will be handled as follows. First, note that
\begin{eqnarray*}
B(r,r') &=& m_{2}(A_{2}) \int_{0} ^{t} \mathrm{d}u \int_{0} ^{t} \mathrm{d}v \int_{A_{2}} \frac{c \mathrm{d}y}{y^{2}} \int_{u}^{t} \mathrm{d}a \int_{u}^{t} \mathrm{d}b \partial _{1} K (a,u) \partial _{1}K(b,u)\\
&&\times \mathbb{E} \left( Q_{r}(a) 1_{C_{r}(a)} (v,y) - Q_{r'}(a) 1_{C_{r'}(a)} (v,y)\right)\left( Q_{r}(b) 1_{C_{r}(b)} (v,y) - Q_{r'}(b) 1_{C_{r'}(b)} (v,y)\right)\\
&=& m_{2}(A_{2}) \int_{0} ^{t} \mathrm{d}u \int_{0} ^{t} \mathrm{d}v \int_{A_{2}} \frac{c \mathrm{d}y}{y^{2}} \int_{u}^{t} \mathrm{d}a \int_{u}^{t} \mathrm{d}b \partial _{1} K (a,u) \partial _{1}K(b,u)\\
&&\times \mathbb{E} \left[ Q_{r}(a) Q_{r}(b) 1_{C_{r}(a)} (v,y) 1_{C_{r}(b)}(v,y) -  Q_{r}(a) Q_{r}(b) 1_{C_{r}(a)} (v,y) 1_{C_{r'}(b)}(v,y)\right. \\
&&\left.
- Q_{r}(a) Q_{r}(b) 1_{C_{r'}(a)} (v,y) 1_{C_{r}(b)}(v,y)+  Q_{r'}(a) Q_{r'}(b) 1_{C_{r'}(a)} (v,y) 1_{C_{r'}(b)}(v,y)\right]
\end{eqnarray*}
where we used again (\ref{mart}). By decomposing
$$C_{r}(a) = C_{r'} (a) \cup \left( C_{r}(a) \setminus C_{r'}(a) \right) $$
we will have
\begin{eqnarray*}
B(r,r') &=& m_{2}(A_{2}) \int_{0} ^{t} \mathrm{d}u \int_{0} ^{t} \mathrm{d}v \int_{A_{2}} \frac{c \mathrm{d}y}{y^{2}} \int_{u}^{t} da \int_{u}^{t} \mathrm{d}b \partial _{1} K (a,u) \partial _{1}K(b,u)\\
&&\times \mathbb{E} \left[ Q_{r'}(a) Q_{r'}(b) -Q_{r}(a) Q_{r} (b) \right] 1_{C_{r'}(a) } (v,y)1_{C_{r'}(b)}(v,y) \\
&&+ m_{2}(A_{2}) \int_{0} ^{t} \mathrm{d}u \int_{0} ^{t} \mathrm{d}v \int_{A_{2}} \frac{c \mathrm{d}y}{y^{2}} \int_{u}^{t} \mathrm{d}a \int_{u}^{t} \mathrm{d}b \partial _{1} K (a,u) \partial _{1}K(b,u)\\
&&\times \mathbb{E}Q_{r}(a) Q_{r} (b) 1_{ C_{r}(a) \setminus C_{r'}(a)}(v,y)  1_{ C_{r}(b) \setminus C_{r'}(b)}(v,y) \\
&:=& B_{1}(r,r')+ B_{2}(r,r').
\end{eqnarray*}
 Obviously,  by bounding the indicator functions by 1,
 \begin{eqnarray*}
 B_{1}(r,r') &\leq & \int_{0} ^{t} \mathrm{d}u \int_{0} ^{t} \mathrm{d}v \int_{A_{2}} \frac{c \mathrm{d}y}{y^{2}} \int_{u}^{t} \mathrm{d}a \int_{u}^{t} \mathrm{d}b \partial _{1} K (a,u) \partial _{1}K(b,u)\\
&&\times \mathbb{E} \left[ Q_{r'}(a) Q_{r'}(b) -Q_{r}(a) Q_{r} (b) \right]
 \end{eqnarray*}
 and therefore it converges to zero as $r,r'\to 0$ by using exactly the same argument as in the case of the term $A(r,r').$

Concerning $B_{2}(r,r')$, since  for $0<r<r'<1$ and $y\in A_{2}$
$$1_{ C_{r}(a) \setminus C_{r'}(a)}(v,y)  \leq 1_{(r,r')\times A_{2}} (v,y)$$
we can bound it in the following way
\begin{eqnarray*}
B_{2}(r,r') &\leq &   m_{2}(A_{2}) \int_{0} ^{t} \mathrm{d}u \int_{0} ^{t} \mathrm{d}v \int_{A_{2}} \frac{c \mathrm{d}y}{y^{2}} \int_{u}^{t} \mathrm{d}a \int_{u}^{t} \mathrm{d}b \partial _{1} K (a,u) \partial _{1}K(b,u)\\
&&\times \mathbb{E}Q_{r}(a) Q_{r} (b)  1_{(r,r')\times A_{2}} (v,y)\\
&=& m_{2}(A_{2}) ^{2} (r'-r) \mathbb{E} \int_{0} ^{t} \mathrm{d}u \int_{u}^{t} \mathrm{d}a \int_{u}^{t} \mathrm{d}b \partial _{1} K (a,u) \partial _{1}K(b,u) \mathbb{E}Q_{r}(a) Q_{r} (b)
\end{eqnarray*}
and by interchanging the order of integration and using
\begin{equation}\label{aux}
\int_{0}^{a\wedge b} \mathrm{d}u\partial _{1}K (a,u) \partial _{1} K(b,u)= c_{H} \vert a-b\vert ^{2H-2}.
\end{equation}
and since
$$\mathbb{E} Q_{r}(a) Q_{r}(b)= e^{-\varphi (2) m_{r}(\vert a-b\vert ) } \leq   e^{c\varphi (2)\log \vert a-b\vert  } = \vert a-b\vert ^{c\varphi (2)}$$
we get
$$B_{2}(r,r') \leq m_{2}(A_{2}) ^{2} (r'-r) \int _{0} ^{t} \int_{0}^{t} \mathrm{d}a \mathrm{d}b  \vert a-b\vert ^{c\varphi (2)+2H-2}$$
we clearly obtain than $B_{2}(r,r')$ goes to zero as $r,r'\to 0$ under condition (\ref{cc}). \qed

\begin{definition}\label{def1}
The process $(Z^{H})_{t>0}$  from Theorem \ref{t1} will be called as  {\it fractional Multifractal Random Walk.  }
\end{definition}

\section{Properties of the MRW}
We will discuss here some immediate properties of the fractional MRW from Definition \ref{def1}. Basically, as we mentioned before, in the disjoint case the fractional MRW $Z^{H}$ has the same properties as in the situation when $Q$ and $B^{H}$ are independent: self-similarity, stationarity of increments and long-range dependence (see \cite{ACCP}). But in the non-disjoint case, we will show that even the $L^{2}$ moment of the fractional MRW does not scale exactly. We provide an exact calculation in order to show this phenomenon. On the other hand, we can control the $L^{p}$ norm of the increment of the process and we find some kind of asymptotic scaling even in the non-disjoint case.

\subsection{Scaling of the second moment}

We first analyze the second moment of the increments of $Z^{H}$.
\vskip0.2cm

{\bf The disjoint case: } Following \cite{ACCP}, we have

\begin{eqnarray*}
\mathbb{E} \left( Z^{H}(t)\right) ^{2}&=&\int_{0}^{t} \mathrm{d}s \left( \int_{s}^{t} \int_{s}^{t} \mathrm{d}a\mathrm{d}b \mathbb{E}(Q_{0}(a) Q_{0}(b))\partial _{1}K (a,u) \partial _{1} K(b,u)\right)\\
&=& \int_{0}^{t} \mathrm{d}s \left( \int_{s}^{t} \int_{s}^{t} \mathrm{d}a\mathrm{d}b  e^{-\rho _{2} m_{0} (\vert a-b\vert )}\partial _{1}K (a,u) \partial _{1} K(b,u)\right)\\
&=& c_{H} \int_{0} ^{t} \int_{0} ^{t} \mathrm{d}b \vert a-b\vert ^{2H-2} \mathbb{E}(Q_{0}(a) Q_{0}(b))
\end{eqnarray*}
where the meaning of the quantity $\mathbb{E}Q_{0}(a) Q_{0}(b)$ is given by (\ref{covQ}) with $r=0$ and where we used the identity (\ref{aux}). Then, for every $h>0$
\begin{eqnarray*}
\mathbb{E} \left( Z^{H}(ht)\right) ^{2}&=& c_{H} \int_{0} ^{ht} \int_{0} ^{ht} \mathrm{d}b \vert a-b\vert ^{2H-2} \mathbb{E}(Q_{0}(a) Q_{0}(b))\\
&=& h^{2H} c_{H} \int_{0} ^{t} \int_{0} ^{t} \mathrm{d}b \vert a-b\vert ^{2H-2} \mathbb{E}(Q_{0}(ha) Q_{0}(hb))\\
&=& h^{2H} \mathbb{E} e^{2 \Omega _{h} } c_{H} \int_{0} ^{t} \int_{0} ^{t} \mathrm{d}b \vert a-b\vert ^{2H-2} \mathbb{E}(Q_{0}(a) Q_{0}(b))\\
&=&h^{2H} \mathbb{E} e^{2 \Omega _{h} } \mathbb{E} \left( Z^{H}(t)\right) ^{2}=h^{2H+ c \varphi (2) }\mathbb{E} \left( Z^{H}(t)\right) ^{2}.
\end{eqnarray*}
where we used the scaling property (\ref{sq}). Actually, it is not difficult to see that for every $p>1$
$$\mathbb{E} \left( Z^{H}(ht)\right) ^{p}=h^{2Hq} \mathbb{E} e^{p\Omega _{t}}=h^{2Hp+ q \varphi (p)} .$$
Moreover, we have the self-similarity $(Z^{H}(at)) _{t\in [0,1] }= ^{(d)} a^{H+\Omega _{a} }(Z^{H}(t))_{t\in [0,1]}$, the stationarity of the increments and the long-range dependence in the sense that $\mathbb{E}X_{k}X_{0}  \sim \tau ^2H k^{2H}$ where $X_{k}= Z^{H}((k+1)\tau )-Z^{H}(k\tau)$ with $k\geq 0$ integer and $\tau >0$.

\vskip0.2cm

{\bf The non-disjoint case: }

The situation is different in the non-disjoint case and we will see that even the second moment of the fractional MRW does not scale.  We can compute exactly the $L^{2}$ norm of $Z^{H}(t)$.

\begin{eqnarray*}
\mathbb{E} \left| Z^{H}(t) \right| ^{2}&=&  \int_{0} ^{t} \mathrm{d}s \left(  \int_{s}^{t}\mathrm{d}a Q_{0}(a) \partial _{1} K(a,s) da \right) ^{2} \\
&&+ \int_{0} ^{t} \mathrm{d}s \int_{A_{2}} \frac{c \mathrm{d}r}{r^{2}} \int_{0} ^{t} \mathrm{d}\alpha \int_{A_{2}} \frac{ c \mathrm{d}\beta } {\beta ^{2}} \\
&&\times \left( D_{\alpha , \beta } \int_{s}^{t} \mathrm{d}b Q_{0}(b) \partial _{1} K(b,s) \mathrm{d}b\right) \left( D_{s,r} \int_{\alpha }^{t} \mathrm{d}a Q_{0}(a) \partial _{1} K(a,\alpha ) \mathrm{d}a\right)\\
&=& \int_{0} ^{t} \mathrm{d}s \left(  \int_{s}^{t}\mathrm{d}a Q_{0}(a) \partial _{1} K(a,s) \mathrm{d}a \right) ^{2} \\
&&+ \int_{0} ^{t} \mathrm{d}s \int_{A_{2}} \frac{c \mathrm{d}r}{r^{2}} \int_{0} ^{t} \mathrm{d}\alpha \int_{A_{2}} \frac{ c \mathrm{d}\beta } {\beta ^{2}}\\
&&\times
\left( \int _{\alpha }^{t} \mathrm{d}aQ_{0} (a) \partial _{1} K (a, \alpha ) 1_{C_{0}(a)}(s,r)  \right) \left( \int_{s}^{t} \mathrm{d}b Q_{0}(b) 1_{C_{0}(b)}(\alpha, \beta)\partial _{1} K(b,s) \mathrm{d}b\right)\\
&=& I(t)+J(t)
\end{eqnarray*}
Consider $h>0$.  Then, applying the above formula to $t=ht$ and making several changes of variables, we will get
\begin{eqnarray*}
\mathbb{E}\left| Z^{H}(ht) \right| ^{2}&=& h^{2H+  c\varphi (2)}  I(t) \\
&&+ h^{2H+1+c\varphi (2)}  \int_{0} ^{t} \mathrm{d}s \int_{A_{2}} \frac{c \mathrm{d}r}{r^{2}} \int_{0} ^{t} d\alpha \int_{A_{2}} \frac{ c \mathrm{d}\beta } {\beta ^{2}}
\\
&&\times \left( \int _{\alpha }^{t} \mathrm{d}aQ_{0} (a) \partial _{1} K (a, \alpha ) 1_{C_{0}(ha)}(hs,r)  \right) \left( \int_{s}^{t} \mathrm{d}b Q_{0}(b) 1_{C_{0}(hb)}(h\alpha, \beta)\partial _{1} K(b,s) \mathrm{d}b\right).
\end{eqnarray*}
The integral with respect to $d\beta $ and $dr$ can be computed explicitly. For example, when $A_{2}= (c, \infty) $, $c<1$, we find
$$\int_{A_{2}} \frac{ c \mathrm{d}\beta } {\beta ^{2}}1_{C_{0}(hb)}(h\alpha, \beta)=c\left[ \frac{1}{c\vee 2h\vert b-\alpha\vert }-1\right] $$
and it is clear that the term $J(t)$ does not scale exactly as $I(t)$.

\subsection{The control of the increments}

Let us estimate the $L^{p}$ norm of the increment $Z^{H} (t) -Z^{H}(s)$ of the fractional  Multifractal Random Walk introduced in Definition 1. Fix $s,t\in [0,1]$ with $t>s$.
We will not insist on the disjoint case because the calculations in \cite{ACCP} still hold, so we will have
\begin{equation*}
\mathbb{E}\left| Z^{H} (t) -Z^{H}(s) \right| ^{p} \sim C_{p} \vert t-s \vert ^{2Hp + c\varphi (p)}
\end{equation*}
with $\varphi$ given by (\ref{fi}).

Let us consider the non-disjoint case. In this case $Z^{H} (t)$, which can formally be written as $\int_{0} ^{t} Q_{0} (y) dB ^{H}(u)$ is an anticipating (Skorohod integral).
We need to use Meyer's inequalities (\ref{meyer})  in order to estimate its $L^{p}$ norm.

Actually

$$Z^{H} (t) -Z^{H}(s)= \int_{0} ^{1} \mathrm{d}W^{(1) }(u) F_{s,t}(u) $$
where we denoted by
\begin{equation*}
F_{s,t} (u)= 1_{(0,t) }(u) \int_{u}^{t} \partial _{1} K (a,u) Q_{0}(a) \mathrm{d}a -1_{(0,s) }(u) \int_{u}^{s} \partial _{1} K (a,u) Q_{0}(a) \mathrm{d}a.
\end{equation*}
By Meyer's inequality (\ref{meyer})
\begin{eqnarray*}
\mathbb{E}\left| Z^{H} (t) -Z^{H}(s) \right| ^{p}&\leq &
\mathbb{E}\int_{0} ^{1}\mathrm{d}u  \int_{A_{2}}\frac{c \mathrm{d}r}{r^{2}} \left|F_{s,t} (u) \right| ^{p} \\
&+& \mathbb{E} \int_{0} ^{1}\mathrm{d}u  \int_{A_{2}}\frac{c \mathrm{d}r}{r^{2}} \int_{0} ^{1}\mathrm{d}\alpha   \int_{A_{2}}\frac{c \mathrm{d}\beta}{\beta ^{2}} \left|  D_{\alpha , \beta} F_{s,t} (u) \right| ^{2} \\
&: =& A(t,s)+ B(t,s).
\end{eqnarray*}
The terms $A$ is exactly the $L^{p}$ norm in the disjoint case, so
$$A(t,s) \leq c_{p} \vert t-s\vert ^{2Hp+c\varphi (p)}.$$
Concerning the summand denoted by $B(t,s)$, we proceed as in the proof of Proposition \ref{p2}: we use first the Malliavin differentiation $D_{\alpha , \beta } Q_{r}(a)= Q_{r}(a) 1_{C_{0}}(\alpha, \beta)$, then we bound the indicator function $1_{C_{0}}(\alpha, \beta)$ by 1, then we integrate $dr$ and $d\beta$ and we obtain
\begin{equation*}
\begin{aligned}
B(t,s) &\leq c_{p} \int_{0} ^{1}\mathrm{d}u \left|F_{s,t} (u) \right| ^{p}\leq c_{p} \vert t-s\vert ^{2Hp+c\varphi (p)}
\end{aligned}
\end{equation*}
because the right hand side is equal, modulo a constant, to $A(t,s)$. Taking into account the above estimates, we conclude that
\begin{equation}
\label{2i}
\mathbb{E}\left| Z^{H} (t) -Z^{H}(s) \right| ^{p}\leq c_{p} \vert t-s\vert ^{2Hp+c\varphi (p)}
\end{equation}
for every $t,s$.

\section{Financial Statistics}

As mentioned in the Introduction, the multifractal random walks appears nowadays as a serious candidate to model the financial time series. In order to compare its behavior with real data, one needs to simulate it. The main difficulty consists in the fact that, in our construction, the variables $Q$ and $B^{H}$ are dependent. From the theoretical point of view, the Malliavin calculus offers  convenient techniques but thinks are complicated in from the practical point of view.

\subsection{Simulation Scheme}
Recall that the fractional MRW is defined by

\begin{equation}
Z(t)=\int_0^tQ(s)\mathrm{d}B^H(s)
=\int_0^t\frac{e^{M(C(t))}}{\mathbb{E}e^{M(C(t))}}\mathrm{d}B^H(s),\label{contt}
\end{equation}

where the measure $M$ and the fBm $B^H$ are dependent.  As throughout our paper, we choose $M$ to be the Gaussian measure.  Consider a random variable $w$ such that $\frac{e^{M(C(t))}}{\mathbb{E}e^{M(C(t))}}=e^{w(t)}$.\\

The first step consists in generating a sequence of random variables $W=(w(t_{i})) _{i}$ such that
 $w(t_{i+1})$ and  $w(t_i)$ are correlated and with Gaussian distribution  $\mathcal{N}_n(\mu, \Sigma)$ where the covariance matrix $\Sigma$ depending on the auto-correlation function   $\gamma(|t - (t+h)|)=Cov(w(t), w(t+h))$. We consider the choice of $\gamma$ from \cite{BMD} with $\gamma(k)$ decreasing with respect to $k$. Then, assuming that the matrix $\Sigma$ is positively definite, there exists a triangular matrix   $C$ such that $C^tC=\Sigma$. In this way, the components of

$$
W=\mu+C^t\mathcal{N}_n(0,\bold{I}),
$$
are correlated and with normal law.

Let us explain the idea to simulate the integral (\ref{contt}). This integral is a divergence integral and in principle its simulation is difficult. But in our particular case, we can use the following approach. Let  $\{0=t_0<t_1<\cdots<t_{n}=t\}$ a partition of the interval $[0,t]$. in the dependent case (meaning when $W $ and $B^{H}$ are independent), (\ref{contt}) can be naturally approximated by

$$
\sum_{i=0}^{n-1}e^{w(t_i)}(B^H(t_{i+1})-B^H(t_{i})).
$$
 Since the simulation of $B^{H}$ is well-known (see e.g. \cite{Bardet} for an explicit algorithm), we can generate the above sum. In the dependent case, using the integration by parts formula $\delta (Fu)= F\delta (u)+ \langle DF, u\rangle$, $\delta $ being the Skorohod integral with respect to $B^{H} $ (see \cite{N}), the same sum can be expressed as
 $$\sum_{i=0}^{n-1}\left( e^{w(t_i)}(B^H(t_{i+1})-B^H(t_{i}))+ \langle De^{w(t_{i})}, 1_{(t_{i},t_{i+1})} \rangle \right)
$$

$$= \sum_{i=0}^{n-1}\left( e^{w(t_i)}(B^H(t_{i+1})-B^H(t_{i}))+  e^{w(t_{i})}\langle Dw(t_{i}) , 1_{(t_{i},t_{i+1})} \rangle \right)$$
where the scalar product is in the Hilbert space associated to the fBm and, since $w(t_{i})$ is Gaussian, $ Dw(t_{i})$ is deterministic and it depends on the set $A_{2}$ (more exactly on the intersection of $A_{2}$ and the interval (0,1) which is in principle very small). Thus, a possibility to simulate (\ref{contt}) is to approximate it by $\sum_{i=0}^{n-1}\left( e^{w(t_i)}(B^H(t_{i+1})-B^H(t_{i}))+a_{i}e^{w(t_{i})}\right)$ with suitable coefficients $a_{i}$.

\subsection{Stylized Facts}

Let $S(t)$ be the price at time $t\ge T$ of an financial asset, $X(t)=\ln S(t)$ the log price and then, log-returns at lag $\tau$ are given by

$$
\delta_\tau X(t)=X(t)-X(t-\tau)=\ln\left(\frac{S(t)}{S(t-\tau)}\right),
$$
and, for sake of simplicity and with no loss of generality, we assume that log-returns are centered, $\mathbb{E}[\delta_\tau X(t)]=0$, thus, we can interpret the volatility as the squared log-returns.\\

The construction of the proposed MRW allows us to verify several properties of stock market fluctuations, what we call stylized facts. The first two that we present are specific to the MRW, long memory in absolute log-return and in volatility; the others  are achieved thanks to our construction, leverage effect and long memory in return.\\

The random walk that we introduce in our work allows the conditional variance and the random noise, which is a fBm with Hurst parameter $H>\frac{1}{2}$, to be dependent. This dependence should involve the so-called leverage effect, that is, the correlation between the log-return at time $t$ and the squared log-return (the volatility) in the future. The case $H>1/2$ implies that the noise $B^H$ exhibits long-range dependence. This should involve a long-range dependence of the fluctuations of the financial time series. These facts (leverage effect and long-range dependence) will be checked on the data. Actually the leverage effect is empirically observed on financial index such as the $S\&P$ 500 (US market), Nikkei 225 (Japan market), FTSE 100 (UK market), CAC 40 (french market) etc. But it does not appear on the Forex market and on the commodities, for example. We also notice the absence of auto-correlation in log-returns if we consider low frequency financial data but this auto-correlation structure clearly appears in the case of square or absolute log-returns. Nevertheless, at high frequency, we observe a dependence relation which is going to be formed even for usual (not squared and not absolute) auto-correlation function. (see e.g. \cite{LF} or \cite{Voit}).\\

\begin{figure}[h!]
	\begin{center}
		\includegraphics[height=6.5cm, width=14cm]{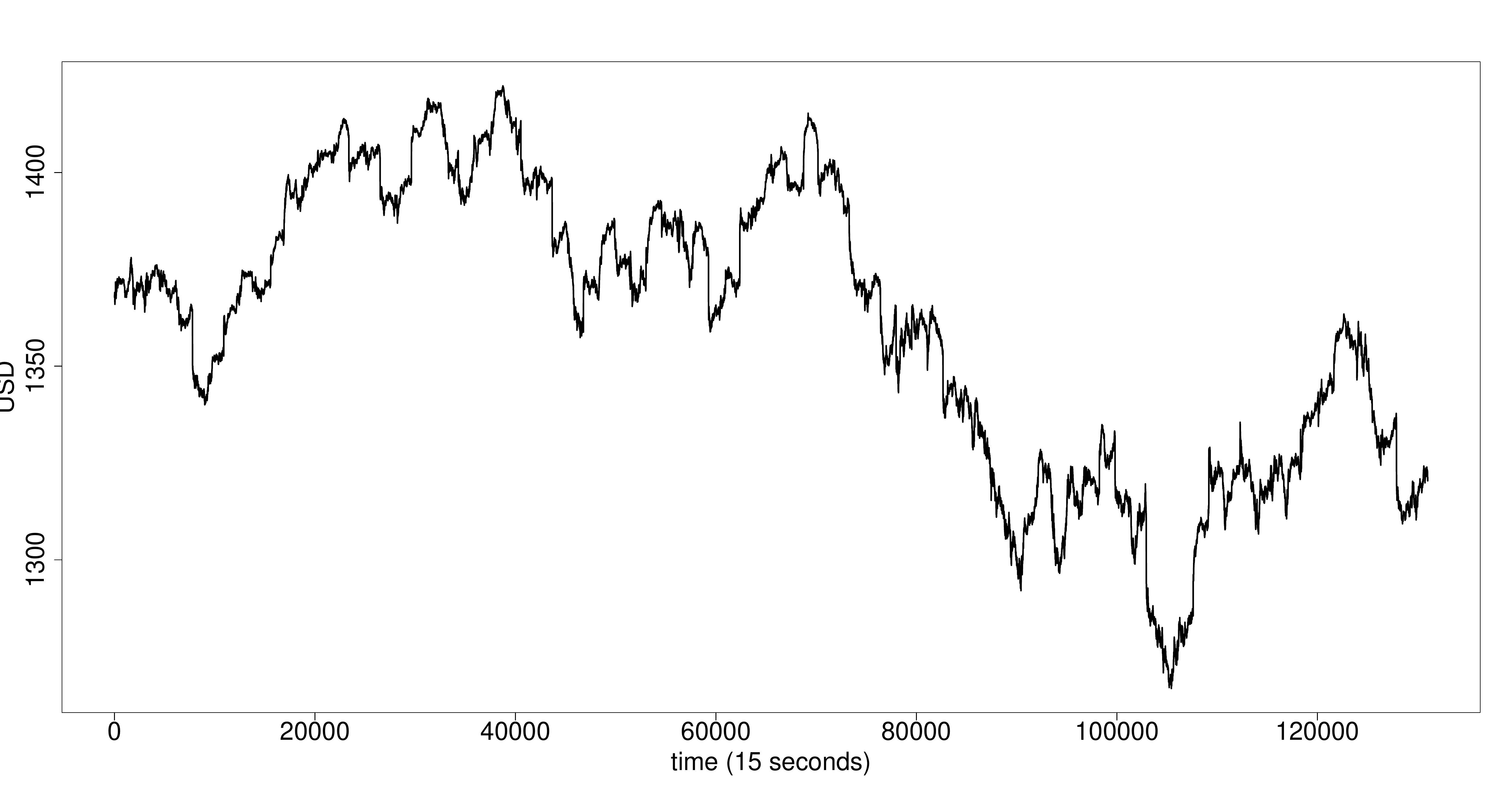}
	\caption{ S$\&$P 500 index , at 15 seconds, from  2012-02-28 to 2012-06-26, 131011 points.}
	\label{plotspy}
	\end{center}
\end{figure}

\begin{figure}[h!]
     \centering \includegraphics[height=6.5cm, width=14cm]{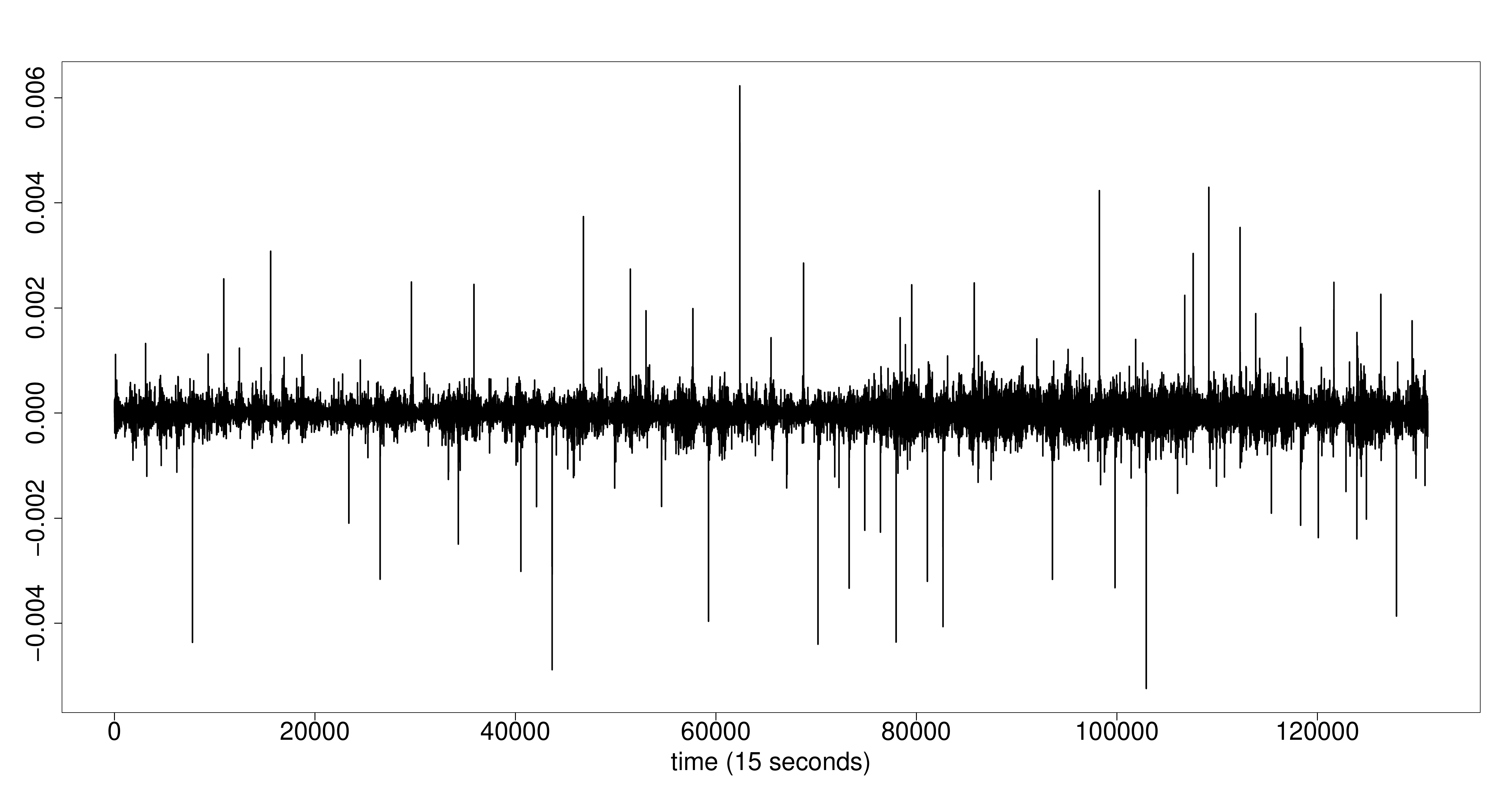}
    \centering \includegraphics[height=6.5cm, width=14cm]{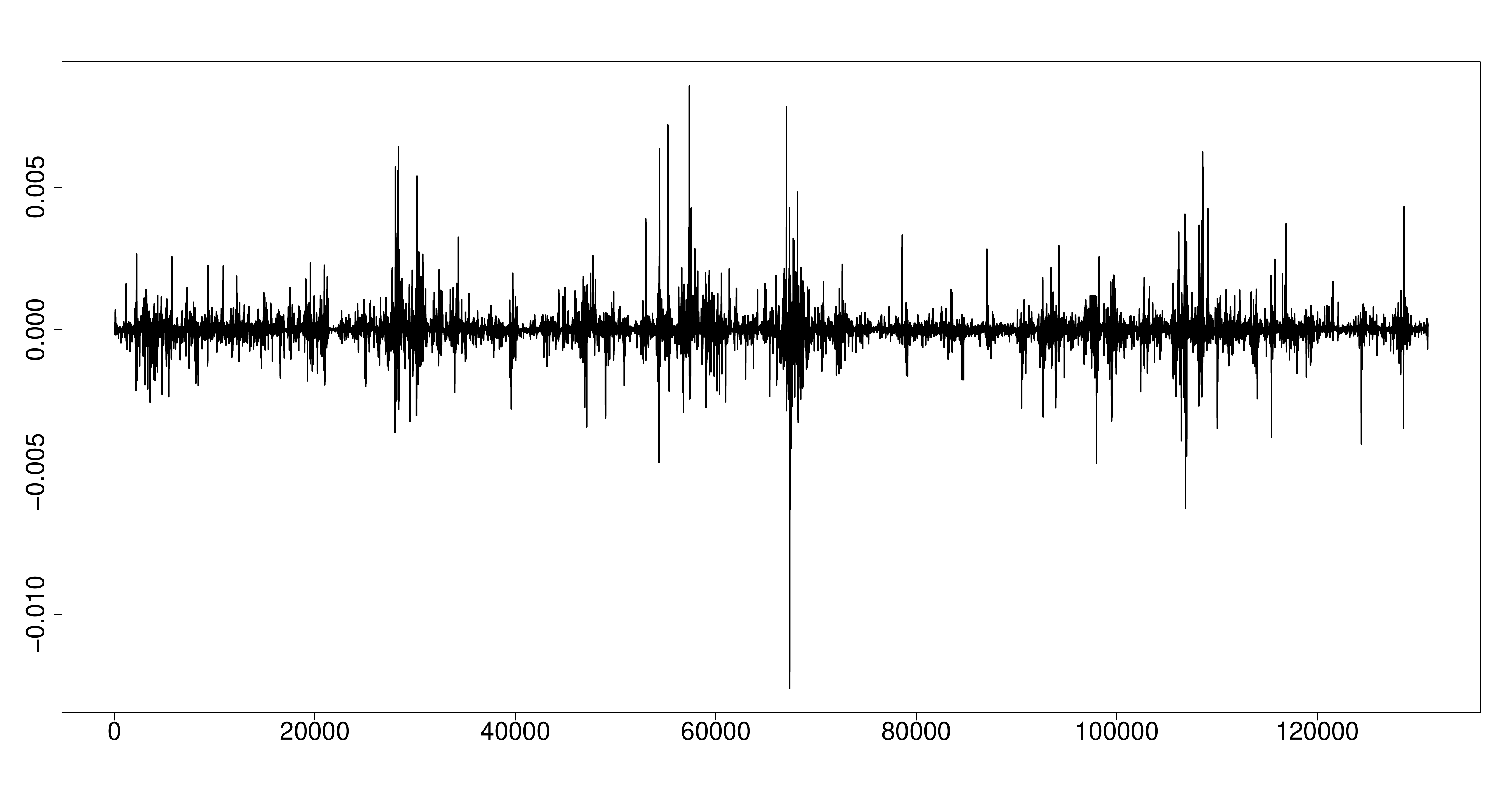}
 \caption{ Top: log-return of S$\&$P 500, at 15 seconds, from  2012-02-28 to 2012-06-26, 131011 points. Down: simulation of the  MRW with $H=0.62$.}
\label{plotreturn}
\end{figure}

The data that we used in order to check the leverage effect and the long-range dependence come from S$\&$P 500 index with a frequency of 15 seconds from 2012-02-28 to 2012-06-26, 131011 points. We compare the characteristics of the S$\&$P 500 index with the simulation of our multifractal random walk. As we mention,  we  compare our simulations with the data from the S$\&$P 500 index. In fact, we only present the comparison for the auto-correlation function since the multifractal structure of the asset return distribution are already well known and well describe by standard multifractal random walk (see e.g. \cite{MDB}). The numerical results are presented in Figures \ref{plotspy} and \ref{plotreturn}. To compare, we also present the simulation of the fractional MRW with $H=0.62$.  This value has been obtained by a classical R/S test.\\

\begin{figure}[h!]
	\begin{center}
		\includegraphics[height=6.5cm, width=14cm]{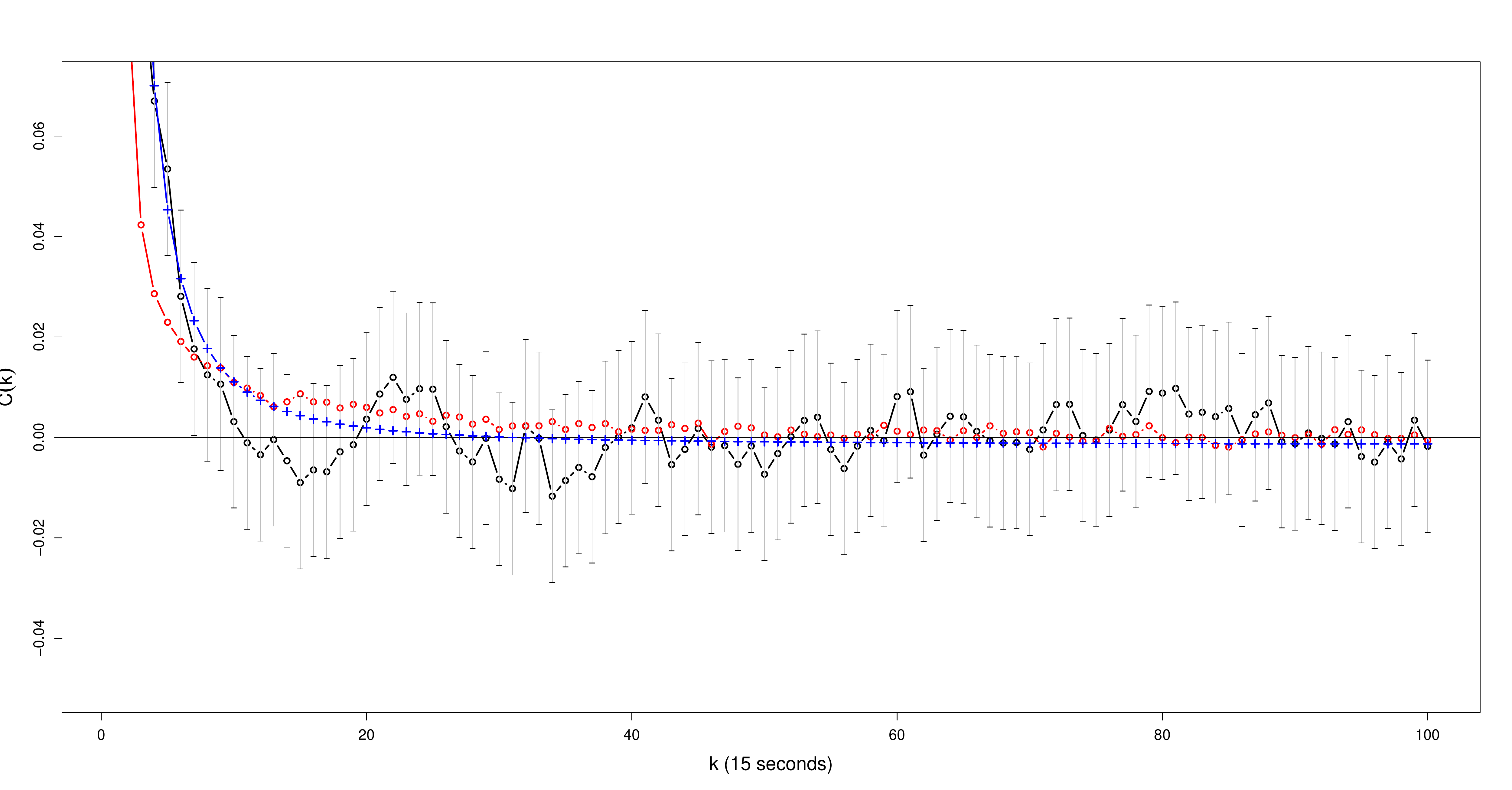}
	\caption{Auto-correlation of log-returns, $k$ from  1 to 100. In black the  S$\&$P 500, in grey the confidence interval, in red the  MRW and in blue to fit in power law ($\alpha=2.1$). Notice that the mean square error between the data and the numerical results is $3.868\times10^{-4}$ and if we take into account the fit, it is $0.404\times10^{-4}$}
	\label{plotacf}
	\end{center}
\end{figure}

Since the Hurst parameter is bigger that $0.5$, we should have a long-range dependence for the log-returns.
The auto-correlation function of the log-returns  is defined by

$$
C(k)=\text{Corr}(\delta_\tau X(t), \delta_\tau X(t+k)).
$$
Usually, we say that the log-return have long-range dependence if  it satisfies,$$
C(k)\sim k^{-\alpha},\:\alpha>0,
$$
We give in Figure  \ref{plotacf}  the empirical auto-correlation function of the log-returns, the best fit obtained in power law  and the simulated auto-correlation function of the fractional MRW with $H=0.62$.\\

We also notice a long-range dependence for the log-returns in the case of S$\&$P 500 index at high frequency, but this less than for the squared or absolute long-returns. We could advance two main ideas to explain this small persistence in high frequency log-returns. First one is that all agents include information, good or bad news, in their investment strategies, but there may have some latencies in the response of investor and cause this auto-correlation. The second one can be cause by huge investors such as mutual fund who want to invest significant sums and place large order. They have to split their large order in a sequence of smaller orders to acquire the number of shares wanted, and then, could cause some persistence. A simple and current strategy to do that is called iceberg order, more sophisticate could be found in the theory of 'optimal trading'.\\

\begin{figure}[h!]
	\begin{center}%
		\includegraphics[height=6.5cm, width=14cm]{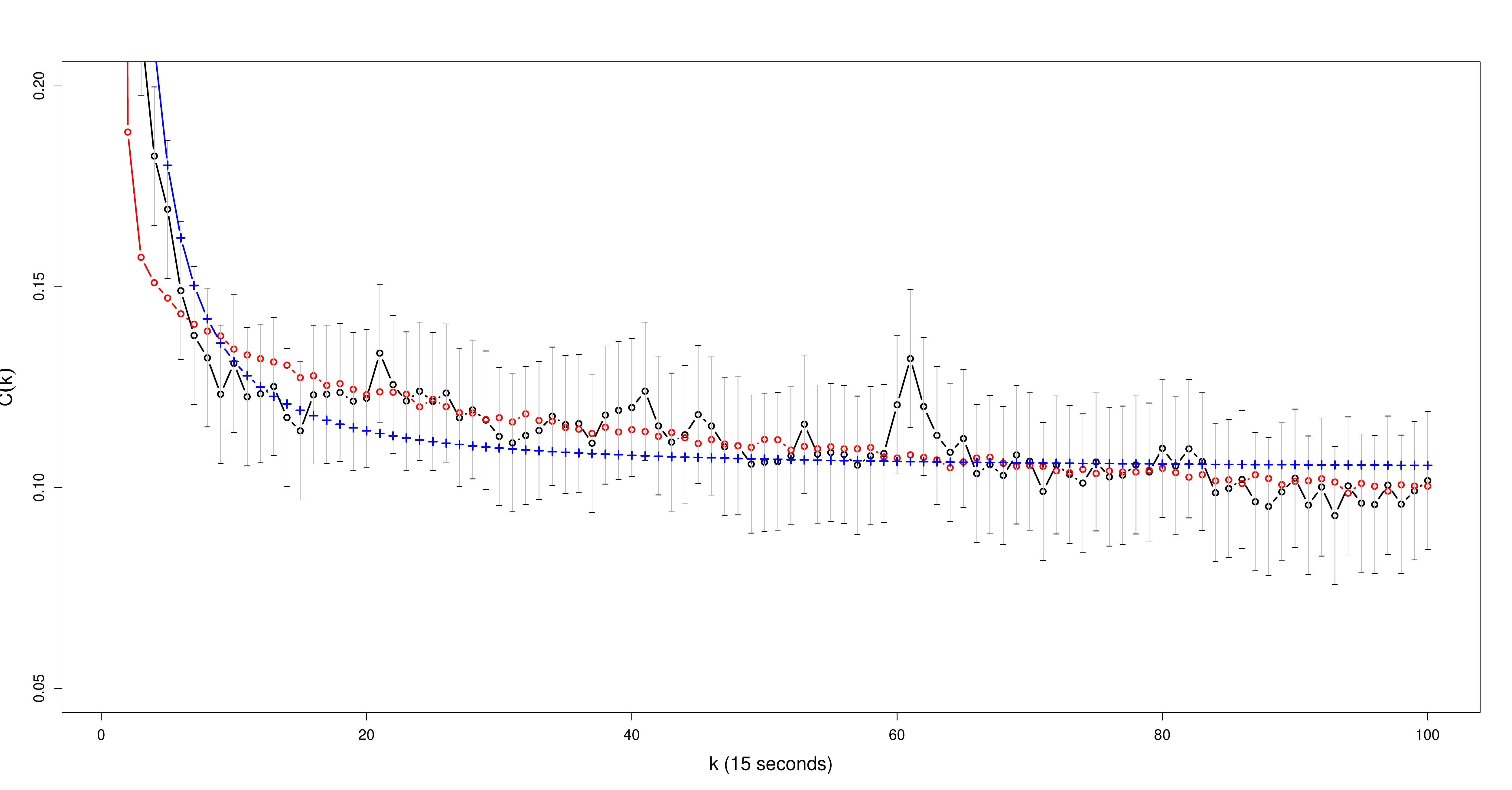}
	\caption{Auto-correlation of the absolute log-returns, $k$ from 1 to 100. In black S$\&$P 500, in grey the confidence interval, in red the  MRW and in blue the fit in power law ($\alpha=1.5$). The mean square error is  $1.720\times10^{-4}$ and with the fit it is $0.223\times10^{-4}$}
	\label{plotacfa}
	\end{center}
\end{figure}

\begin{figure}[h!]
	\begin{center}
		\includegraphics[height=6.5cm, width=14cm]{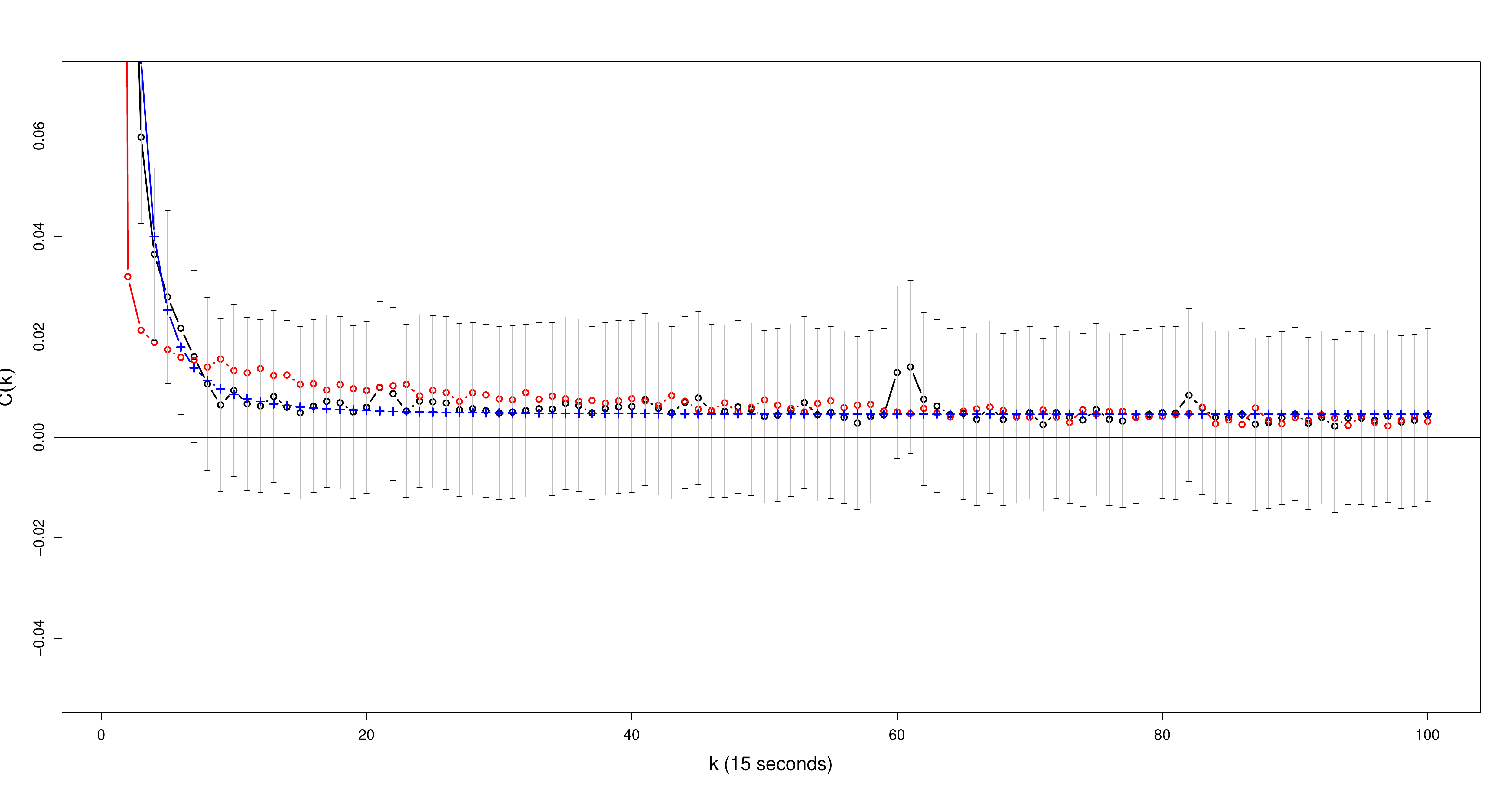}
	\caption{Auto-correlation of the squared log-returns, $k$ from 1 to 100. In black S$\&$P 500, in grey the confidence interval, in red the  MRW and in blue the fit in power law  ($\alpha=2.4$). The mean square error is  $2.834\times10^{-4}$ and with the fit $1.209\times10^{-2}$}
	\label{plotacfs}
	\end{center}
\end{figure}

 The persistence of dependence for the log-returns it not easily observed, but the presence of the cluster of log-returns indicates basically that strong variations are followed by other strong variations. Then, we should observe a persistence in absolute and squares log-returns. As before, we present in Figure \ref{plotacfs} and Figure \ref{plotacfa} the empirical results for $S\& P$ 500 with its fit in power law and the simulation of the MRW ($H$ still equal to  0.62).\\

We now present the second result that we obtain via Malliavin calculus construction. We have assumed that the measure $M$ and the integrator $B^H$ was dependent, hence, the MRW should verify the correlation between return at time $t$ and squared return (volatility) around $t$, say $t+k\leq t'\leq t+k$, this quantity correspond to the leverage effect,

\begin{equation}
\mathcal{L}(t,t')=\mathbb{E}\left[\delta_\tau X(t)(\delta_\tau X(t'))^2\right], \quad t+k\leq t'\leq t+k,
\end{equation}
This quantity is usually normalize by $\mathbb{E}[(\delta_\tau X(t))^2]^2$, we use this normalized formulation, see \cite{BMP} for details. In the literature, one can usually find studies of the stylized facts at low frequency, based on financial indices or on the average of financial indices (see \cite{BLT} for a very interesting empirical study), that's why we have chosen to work with the S$\&$P 500. We will follow the same approach, but we also present results for high frequency.

\begin{figure}[h!]
	\begin{center}
		\includegraphics[height=6.5cm, width=14cm]{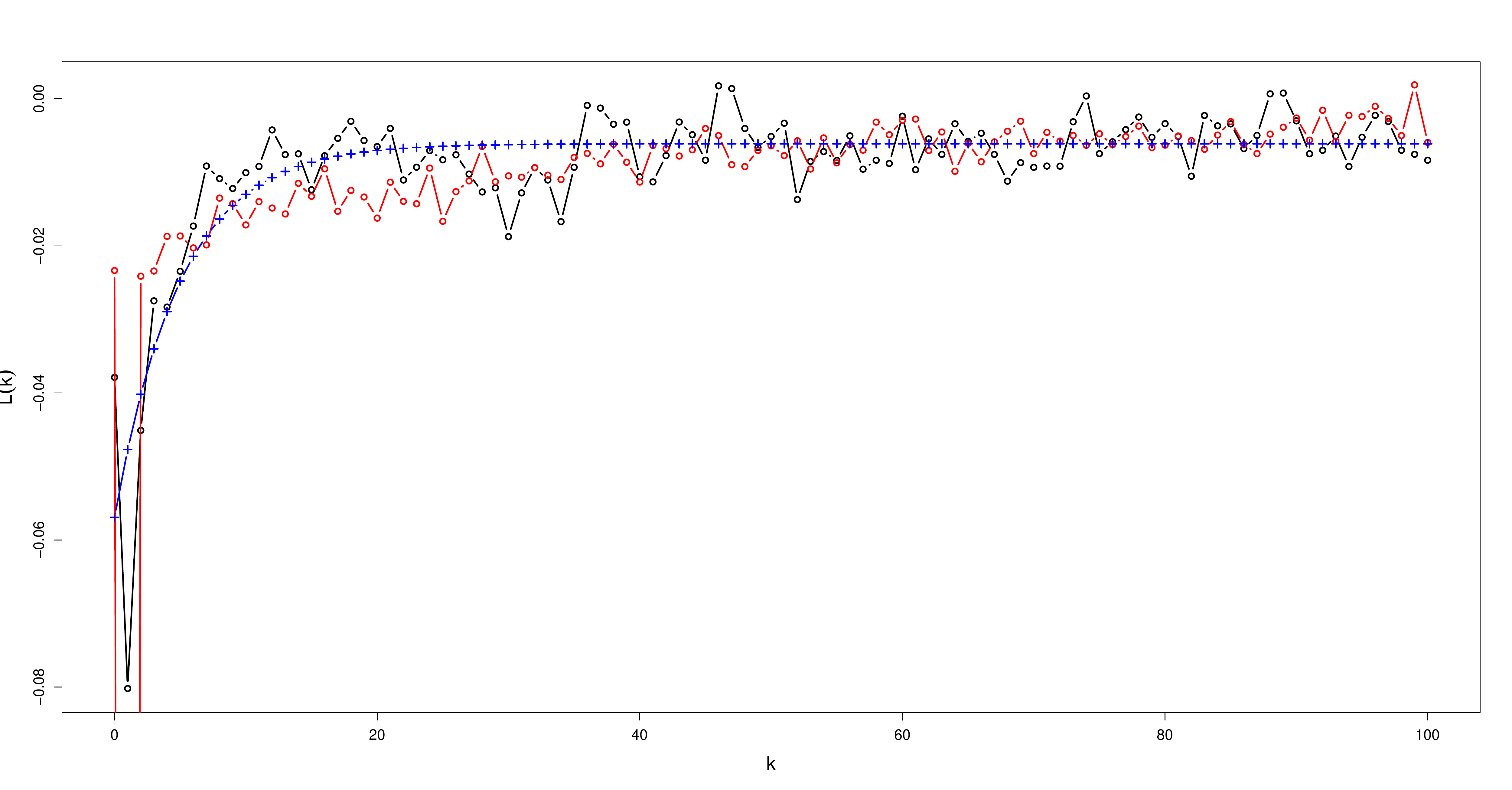}
		\includegraphics[height=6.5cm, width=14cm]{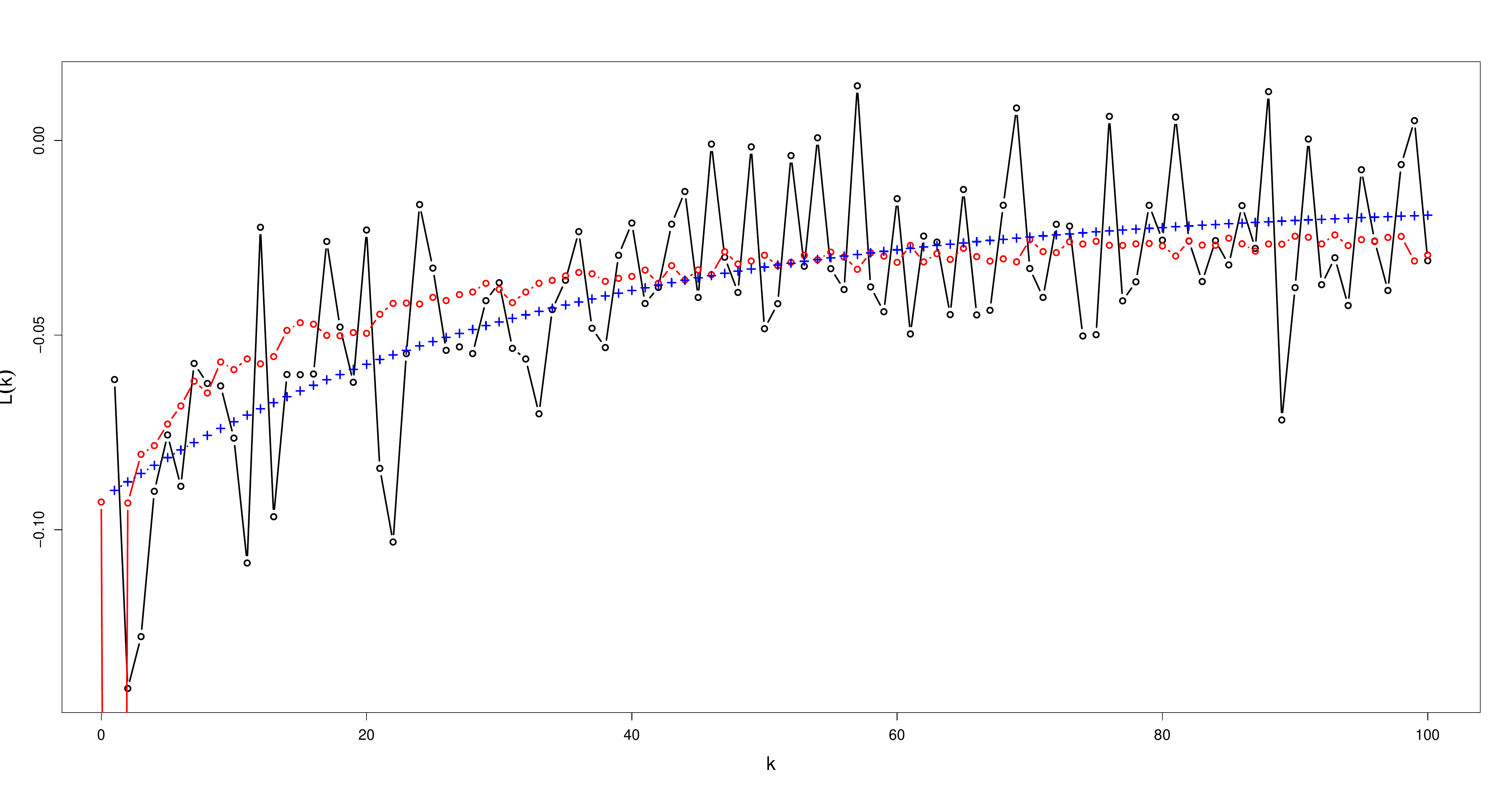}
	\caption{First, the empirical leverage effect on  S$\&$P 500 at  5 seconds in black, best fit in blue, with  $\alpha=0.2$, in red. Second, the leverage effect  for daily data on  S$\&$P 500, from  1985-02-10 to 2012-06-26, the fit is obtained for $\alpha=0.03$, here $H=0.5$.}
	\label{plotlev}
	\end{center}
\end{figure}

The empirical study of the leverage effect is given in Figure \ref{plotlev}. To compare, we also plot the leverage effect for  the fractional MRW in Figure \ref{plotlev} at low frequency.

 \end{document}